\documentclass[11pt]{article}

\usepackage{amsmath}
\usepackage{amssymb}
\usepackage{eucal}

\begin{document}

\baselineskip 15pt

\title{On $\sigma$-semipermutable subgroups of finite groups\thanks{Research is supported by
a NNSF grant of China (Grant \# 11371335) and Wu Wen-Tsun Key Laboratory of Mathematics of Chinese Academy of Sciences.}}

\author{Wenbin Guo\\
{\small Department of Mathematics, University of Science and
Technology of China,}\\ {\small Hefei 230026, P. R. China}\\
{\small E-mail:
wbguo@ustc.edu.cn}\\ \\
{ Alexander  N. Skiba
}\\
{\small Department of Mathematics,  Francisk Skorina Gomel State University,}\\
{\small Gomel 246019, Belarus}\\
{\small E-mail: alexander.skiba49@gmail.com}}

\date{}
\maketitle

\begin{abstract} Let   $\sigma =\{\sigma_{i} | i\in I\}$ be some
 partition of the set of all primes $\Bbb{P}$, $G$ a  finite group and $\sigma (G)
=\{\sigma_{i} |\sigma_{i}\cap \pi (G)\ne  \emptyset \}$. A set
 ${\cal H}$ of subgroups of $G$  is said to be a  \emph{complete
Hall $\sigma $-set} of $G$  if every member $\ne 1$ of  ${\cal H}$
 is a Hall $\sigma _{i}$-subgroup of $G$ for some $\sigma _{i}\in \sigma $ and
   ${\cal H}$ contains exact one Hall  $\sigma _{i}$-subgroup of $G$ for every
 $\sigma _{i}\in  \sigma (G)$.
A subgroup $H$ of $G$ is said to be:  \emph{${\sigma}$-semipermutable in $G$
 with respect to ${\cal H}$} if $HH_{i}^{x}=H_{i}^{x}H$
for all $x\in G$ and all $H_i\in {\cal H}$ such that
 $(|H|, |H_{i}|)=1$; \emph{${\sigma}$-semipermutable in $G$}
 if $H$ is ${\sigma}$-semipermutable in $G$ with respect to some complete Hall
 $\sigma $-set of $G$.

We study the structure of $G$ being based on the
 assumption that some subgroups of $G$ are ${\sigma}$-semipermutable in $G$.

\end{abstract}


\let\thefootnote\thefootnoteorig

\section{Introduction}

Throughout this paper, all groups are finite and $G$ always denotes a
 finite group.  Moreover,  $\mathbb{P}$ is the set of all    primes,
$p\in \pi \subseteq  \Bbb{P}$ and  $\pi' =  \Bbb{P} \setminus \pi$. If $n$ is an integer,
 the symbol $\pi (n)$ denotes the set of all primes dividing $n$; as usual,  $\pi (G)=\pi (|G|)$, the set of all primes dividing the order of $G$.

In what follows, $\sigma =\{\sigma_{i} | i\in I\subseteq \mathbb{N} \}$ is
  some partition of $\mathbb{P}$, that is, $\Bbb{P}=\cup_{i\in I} \sigma_{i}$ and $\sigma_{i}\cap
\sigma_{j}= \emptyset  $ for all $i\ne j$.  Let
 $\sigma (G) =\{\sigma_{i} |\sigma_{i}\cap \pi (G)\ne  \emptyset  \}$.

In the mathematical practice, we often  deal with the
 following two special partitions of  $\Bbb{P}$:
 $\sigma =\{\{2\}, \{3\},  \ldots \}$ and  $\sigma  =\{\pi, \pi'\}$
  (in particular, $\sigma =\{\{p\},
\{p\}'\}$, where $p$ is a prime).

A  set  ${\cal H}$ of subgroups of $G$ is a
 \emph{complete Hall $\sigma $-set} of $G$  \cite{1, 111} if
 every member $\ne 1$ of  ${\cal H}$ is a Hall $\sigma _{i}$-subgroup of $G$
 for some $\sigma _{i} \in \sigma$ and ${\cal H}$ contains exact one Hall
 $\sigma _{i}$-subgroup of $G$ for every  $\sigma _{i}\in  \sigma (G)$.

Subgroups $A$ and $B$ of $G$ are called \emph{permutable} if $AB=BA$.
 In this case they also say that $A$ \emph{permutes} with $B$.

{\bf Definition 1.1.}  Suppose that  $G$ possesses  a complete Hall $\sigma$-set $
  {\cal H}=\{H_{1}, \ldots , H_{t} \}$.  A subgroup $H$ of $G$ is said to be:
  \emph{${\sigma}$-semipermutable in $G$ with respect to ${\cal H}$} if
 $HH_{i}^{x}=H_{i}^{x}H$ for all $x\in G$ and all $i$ such that
$(|H|, |H_{i}|)=1$; \emph{${\sigma}$-semipermutable in $G$}
 if $H$ is ${\sigma}$-semipermutable in $G$ with respect to some complete Hall
 $\sigma$-set of $G$.

Many known results deal with two  special cases of the
  $\sigma$-semipermutability condition: when
 $\sigma =\{\{2\}, \{3\},  \ldots \}$ and   $\sigma =\{\pi, \pi'\}$.

Consider some typical examples.

{\bf Example 1.2.}   A subgroup
 $H$ of $G$ is said to be \emph{$S$-semipermutable} in $G$ if $H$  permutes
 with all Sylow subgroups $P$ of $G$  satisfying  $(|H|, |P|)=1$.  Thus
$H$ is $S$-semipermutable in $G$ if and only if it is
${\sigma}$-semipermutable in $G$ where $\sigma =\{\{2\}, \{3\},  \ldots
\}$.

The $S$-semipermutability condition  can be found in many known  results
 (see for example Section 3 in \cite[VI]{hupp}, Chapter 3 in \cite{GuoII}
 and also the recent papers   \cite{III, Isaacs, Isber}).

Before continuing,  let's make the following remark.

{\bf Remarks 1.3.}
 Let $G=AB$  by a product of subgroups $A$ and $B$ and
$K\leq B$. Suppose that  $A$ permutes with $K^{b}$ for all $b\in B$. Then:

(i)  For any $x=ab$, where
$a\in A$ and $b\in B$, we have $AK^{x}=Aa(K^{b})a^{-1}=a(K^{b})a^{-1}A=K^{x}A$ and hence
$A$ permutes with all conjugates of $K$.

(ii) $A^{x}K=KA^{x}$ for all $x\in G$.  Indeed,
$(A^{x}K)^{x^{-1}}=AK^{x^{-1}}=K^{x^{-1}}A$ by Part (i), so
$(AK^{x^{-1}})^{x}= A^{x}K=KA^{x}$.

{\bf Example 1.4.}   A subgroup $H$ of $G$ is said to be \emph{$SS$-quasinormal}
 if $G$ has a subgroup $T$ such that $HT=G$ and $H$ permutes with all Sylow
 subgroups of $T$.  If  $P$ is a Sylow subgroups  of $T$  satisfying  $(|H|,
|P|)=1$, then $P$  is a Sylow subgroups  of $G$ and so $H$   is
 ${\sigma}$-semipermutable in $G$, where $\sigma =\{\{2\}, \{3\},  \ldots
\}$, by
 Example 1.2 and Remark 1.3(i).  Various applications of $SS$-quasinormal subgroups can
be found  in \cite{shirong, LiSK, KKK} and in many other papers.

{\bf  Example 1.5.}   In \cite{h2},  Huppert proved that if a Sylow $p$-subgroup $P$ of $G$
 of order
 $|P| > p$ has a complement $T$ in $G$ and $T$ permutes with all maximal subgroups of $P$, then
$G$ is $p$-soluble.  In view of Remark 1.3 the condition
 "$T$ permutes with all maximal subgroups of $P$" is equivalent to the
condition "all maximal subgroups of $P$ are ${\sigma}$-semipermutable in $G$ with
 respect to $\{P, T\}$",
 where $\sigma =\{\{p\}, \{p\}'\}$.  The result of Huppert was developed in the papers
 \cite{Ser, Bor}, where instead of maximal subgroups were considered the subgroups of $P$
 of fixed order $p^{k}$.

Further,  the results in \cite{h2, Ser, Bor} were generalized
in  \cite{GuoSS, GuoS163}, where instead of a Sylow $p$-subgroup of $G$  was considered a Hall
subgroup of $G$ (see Section 4 below).

Finally, note  that all the above-mentioned results deal with two special
cases: a "binary" case, when $\sigma =\{\pi, \pi'\}$, and
an  "n-ary" case,  when  $\sigma =\{\{2\}, \{3\},  \ldots
\}$.

In this paper,  we consider the ${\sigma}$-semipermutability condition for  arbitrary
 partition ${\sigma}$ of  $\Bbb{P}$.

In fact, our main results are the following two observations.

{\bf Theorem A. } {\sl  Let  $P$ be a Sylow $p$-subgroup of $G$. Suppose that
 $G$ has a complete Hall $\sigma$-set ${\cal H}=\{H_{1}, \ldots
, H_{t} \}$  such that $H_{1}$  is  $p$-supersoluble  of order divisible by $p$.
 Suppose also that there is a natural number
$k$ such that $p^{k} < |P|$ and every subgroup of $P$ of order $p^{k}$ and every
 cyclic subgroup of $P$ of order 4 (if $p^k=2$ and $P$ is non-abelian) are $\sigma $-semipermutable in $G$  with respect to ${\cal H}$.
 Then $G$ is $p$-supersoluble. }

{\bf Theorem B.} {\sl  Let $X\leq E$ be  normal subgroups of $G$.
 Suppose that $G$ has a complete Hall $\sigma$-set ${\cal H}$ such that
every member of  ${\cal H}$ is  supersoluble.
 Suppose  also that for every non-cyclic Sylow subgroup $P$ of $X$ there is a natural number
 $k=k(P)$ such that $p^{k} < |P|$ and every subgroup of $P$ of order $p^{k}$ and every
 cyclic subgroup of $P$ of order 4 (if $p^k=2$ and $P$ is  non-abelian) are  $\sigma $-semipermutable in $G$  with respect to ${\cal H}$. If
 $X=E$ or $X=F^{*}(E)$, then every   chief factor of $G$ below $E$ is cyclic.  }

In this theorem $F^{*}(E)$ denotes the generalized Fitting subgroup
 of $E$, that is, the product of all normal quasinilpotent subgroups of $E$.

We prove Theorems A and B in Section 3.   In Section 4 we discuss some
applications of these two results.

All unexplained notation and terminology are standard. The reader is referred to
 \cite{DH},  \cite{Bal-Ez}, \cite{prod} or \cite{GuoII} if necessary.

\section{Base lemmas}

Suppose that $G$ has a  complete Hall $\sigma$-set ${\cal H}=\{H_{1}, \ldots , H_{t} \}$.
 For any subgroup $H$ of $G$ we write $H\cap {\cal H}$ to denote the set $\{H\cap
H_{1}, \ldots , H\cap H_{t} \}$. If $H\cap {\cal H}$ is
 a complete Hall $\sigma$-set of $H$, then we say that ${\cal H}$ \emph{reduces into} $H$.

{\bf Lemma  2.1.} {\sl  Suppose that $G$ has a  complete Hall
 $\sigma$-set ${\cal H}=\{H_{1}, \ldots , H_{t} \}$ such that a subgroup
  $H$ of $G$  is   ${\sigma}$-semipermutable  with respect to ${\cal H}$.
  Let  $R$ be a normal subgroup of $G$ and $H\leq L\leq G$. Then:}

(1) {\sl ${\cal H}_{0}=\{H_{1}R/R, \ldots , H_{t}R/R \}$ is a complete Hall
 $\sigma$-set of $G/R$. Moreover, if   for every prime $p$ dividing $|H|$ and for
 a  Sylow $p$-subgroup $H_{p}$
of $H$ we have $H_{p}\nleq R$, then $HR/R$ is ${\sigma}$-semipermutable in $G/N$
  with respect to ${\cal H}_{0}$.}

(2) {\sl If  ${\cal H}$ reduces into $L$, then  $H$ is ${\sigma}$-semipermutable
 in $L$ with respect to $L\cap {\cal H}$. In particular, if $L$ is normal in $G$,
 then $H$ is ${\sigma}$-semipermutable  in $L$ with respect to $L\cap {\cal H}$. }

(3) {\sl If $L\leq H_{i}$, for some $i$, then ${\cal H}$ reduces into $LR$.}

(4) {\sl If $H\leq H_{i}$, for some $i$,,  then $H$ is ${\sigma}$-semipermutable in $HR$.}

(5) {\sl If   $H$ is a $p$-group, where  $p\in  \pi (H_{i})\subseteq
\sigma _{i}$
 and $R$ is
a $\sigma _{i}$-group, then $|G:N_{G}(H\cap R)|$ is a $\sigma _{i}$-number.}

{\bf Proof.}  Without loss of the  generality  we can assume that $H_{i}$
is a $\sigma_{i}$-group for all $i=1, \ldots , t$.

(1) It is clear that ${\cal H}_{0}=\{H_{1}R/R, \ldots , H_{t}R/R \}$ is
  a  complete Hall $\sigma$-set of  $G/R$.
Let  $i\in \{1, \ldots , t\}$  such that $(|HR/R|, |H_{i}R/R|)=1$. Let $p\in \pi (H)$ and
 $H_{p}$ a Sylow $p$-subgroup of $H$. Assume that $p$ divides $|H_{i}|$. Then $H_{i}$
 contains a Sylow $p$-subgroup  of $G$ since it is a Hall subgroup of $G$ and so
 $H_{p}\leq R$,  contrary to the hypothesis. Hence $(|H|, |H_{i}|)=1$. By hypothesis,
 $HH_{i}^{x}=H_{i}^{x}H$ for all $x\in G$. Then
 $$(HR/R)(H_{i}R/R)^{xR}=HH_{i}^{x}R/R$$ $$=H_{i}^{x}HR/R=
(H_{i}R/R)^{xR}(HR/R),$$ so  $HR/R$ is  ${\sigma}$-semipermutable in
  $G/R$
  with respect to ${\cal H}_{0}$.

(2) Let $L_{i}=H_{i}\cap L$ for all $i=1, \ldots , t$ and ${\cal L}=\{L_1, \ldots, L_t\}$.
  By hypothesis, $\cal L$ is a complete $\sigma$-Hall set of $L$.   Let  $i\in \{1, \ldots ,
 t\}$  such that $(|H|, |L_{i}|)=1$ and let  $a\in L$. Then $(|H|, |H_{i}|)=1$. Hence, by
 hypothesis,  $HH_{i}^{a}=H_{i}^{a}H$ for all $a\in L$, so $L\cap
HH_{i}^{a}=H(L\cap H_{i}^{a})=
 H(L\cap
 H_{i})^{a}= HL_{i}^{a}=L_{i}^{a}H$. This shows that $H$ is ${\sigma}$-semipermutable in
 $L$ with respect to $L\cap {\cal H}$.

(3) Since $H_{i}\cap R$ is a Hall  $\sigma_{i}$-subgroup of $R$  and $H_{i}\cap LR=
 L(H_{i}\cap R)$, we have $|LR:H_{i}\cap LR|=|R:H_{i}\cap R|$. Hence $H_{i}\cap LR$
 is a Hall $\sigma_{i}$-subgroup of $LR$. It is clear also that $H_{j}\cap
LR=H_{j}\cap R$ is a Hall $\sigma_{j}$-subgroup of $LR$
 for all $j\ne i$. Hence ${\cal H}$ reduces into $LR$.

(4) This follows from Parts (2) and (3).

(5) For any $j\ne i$,  $H_{j}H=HH_{j}$ is a subgroup of $G$  and
 $HH_{j}\cap R=(H\cap R)(H_{j}\cap R)=H\cap R$, so $H_{j}\leq N_{G}(H\cap 
R)$. Hence
 $|G:N_{G}(H\cap R)|$ is a $\sigma_{i}$-number.

{\bf Lemma 2.2} (See  Kegel \cite{KegI}). {\sl Let $A$ and $B$ be subgroups of $G$ such
 that  $G\ne AB$ and $AB^{x}=B^{x}A$, for all $x\in G$. Then $G$ has a proper normal
 subgroup $N$ such that either $A\leq N$ or $B\leq N$.}

{\bf Lemma 2.3.} {\sl  Let $P$ be a Sylow $p$-subgroup of $G$ and
 ${\cal H}=\{H_{1}, \ldots , H_{t} \}$   a complete Hall $\sigma$-set of $G$ such that
  $p\in \pi (H_{1})$. Suppose that for any $x\in G$, $P^{x}H_{i}$ is a $p$-soluble
 subgroup of $G$ for all $i =2, \ldots , t$.  Then $G$ is $p$-soluble. }

{\bf Proof.}  Assume that this is false and let $G$ be a  counterexample of minimal order
. First note  that the hypothesis  holds for every normal subgroup $R$ of $G$.
Therefore every proper normal subgroup of $G$ is $p$-soluble by the choice of $G$. Moreover,
 the choice of $G$ and the hypothesis imply that $PH_{i}\ne G$ for all $i =2, \ldots , t$.  By Lemma 2.2, we have either $P^{G}\ne G$ or $(H_{2})^{G}\ne G$. Hence $G$ has a proper non-identity normal subgroup $R$. But then $R$ is $p$-soluble. On the other hand, the hypothesis holds for $G/R$, so $G/R$ is also $p$-soluble by the choice of $G$. This implies that $G$ is $p$-soluble.

A group $G$ is said to be  \emph{strictly $p$-closed}
 \cite[p.5]{We} whenever $G_{p}$, a Sylow $p$-subgroup of $G$, is
 normal in $G$ with $G/G_{p}$ abelian of exponent dividing $p-1$. A normal subgroup $H$ of $G$
  is called
\emph{hypercyclically embedded} in $G$ if every chief factor of $G$ below
$H$ is cyclic.

{\bf Lemma 2.4}  {\sl  A normal $p$-subgroup $P$ of $G$ is hypercyclically embedded in $G$ if and only if $G/C_{G}(P)$ is strictly $p$-closed.}

{\bf Proof.} If $P$ is hypercyclically embedded in $G$, then for any chief factor $H/K$ of $G$ below $P$, $G/C_{G}(H/K)$ is abelian of exponent dividing $p-1$. Hence $G/C$, where $C$ the intersection the centralizers of all such factors, is also an abelian group of exponent dividing $p-1$.  On the other hand, $C/C_{G}(P)$ is a $p$-group by \cite[Ch.5, Corollary 3.3]{Gor}. Hence $G/C_{G}(P)$ is strictly $p$-closed.

Now assume that $G/C_{G}(P)$ is strictly $p$-closed and let $H/K$ be any chief
 factor below $P$. Since $C_{G}(P)\leq C_{G}(H/K)$, $G/C_{G}(H/K)$ is strictly
 $p$-closed. But since
$O_{p}(G/C_{G}(H/K))=1$ \cite[Ch.A, Lemma 13.6]{DH}, $G/C_{G}(H/K)$ is abelian
 of exponent dividing $p-1$. It follows from \cite[Ch.1, Theorem 1.4]{We} that $|H/K|=p$.
 Thus $P$ is hypercyclically embedded in $G$.

Let $P$ be a $p$-group.  If $P$  is not a non-abelian 2-group, then we use $\Omega (P)$ to denote  the subgroup $\Omega _{1}(P)$. Otherwise, $\Omega (P)= \Omega _{2}(P)$.

{\bf Lemma 2.5} (See \cite[Lemma 2.12]{2014}). {\sl Let  $P$ be
 a normal $p$-subgroup of $G$ and $D={\Omega}(C)$, where $C$ is  a Thompson critical subgroup of $P$. If either  $P/\Phi (P)$ is hypercyclically embedded in $G/\Phi (P)$   or $D$ is hypercyclically embedded in $G$, then $P$ is also hypercyclically embedded in $G$. }

{\bf Lemma 2.6.} {\sl  Let $C$ be a Thompson  critical subgroup of a $p$-group $P$. Then the group $D={\Omega}(C)$ is of exponent $p$ if $p$ is odd prime or exponent 4 if $P$ is non-abelian 2-group. Moreover, every   non-trivial $p'$-automorphism of $P$ induces a non-trivial automorphism  of $D$.}

{\bf Proof.} The first assertion follows from \cite[Ch.5, Theorem 3.11]{Gor} and \cite[Lemma 2.11]{2014}. The second one directly follows from \cite[Ch. 5, Theorem 3.11]{Gor}.

{\bf Lemma  2.7.}  {\sl Let $E$ be a normal subgroup of $G$ and $P$ a Sylow $p$-subgroup of $E$ such that $(p-1, |G|)=1$.  If either   $P$ is cyclic or $G$ is
 $p$-supersoluble, then $E$ is $p$-nilpotent and $E/O_{p'}(E)\leq Z_{\infty}(G/O_{p'}(E))$. }

{\bf Proof.} First note that in view of \cite[Ch.IV, Theorem 5.4]{hupp} and the condition $(p-1, |G|)=1$, $E$ is $p$-nilpotent.  Let $H/K$ be any chief factor of $G$ such that
 $O_{p'}(E)\leq K < H\leq E$. Then $|H/K|=p$, so $G/C_{G}(H/K)$ divides $p-1$.   But by hypothesis, $(p-1, |G|)=1$,  so $C_{G}(H/K)=G.$ Thus  $E/O_{p'}(E)\leq Z_{\infty}(G/O_{p'}(E))$.

The following lemma is well-known (see, for example, \cite[Lemma 2.1.6]{prod}).

{\bf Lemma 2.8.}  {\sl  If $G$ is $p$-supersoluble and $O_{p'}(G)=1$, then   $p$ is the largest prime dividing $|G|$, $G$ is supersoluble and $F(G)=O_{p}(G)$ is a  Sylow $p$-subgroup of $G$.}

{\bf Lemma 2.9} (See  \cite{knyag}). {\sl Let $H$, $K$  and $N$ be
 pairwise permutable subgroups of $G$ and $H$ is a Hall subgroup of $G$,  then $N\cap HK=(N\cap H)(N\cap K).$}

The following fact is also well-known (see for example \cite[Ch.1, Lemma 5.35(6)]{GuoII}).

{\bf Lemma 2.10} {\sl If $H$ is a subnormal $\pi $-subgroup of $G$, then $H\leq O_{\pi}(G).$}

{\bf Lemma 2.11} (See \cite[Theorem C]{Sk1}). {\sl Let $E$ be a normal subgroup of $G$. If $F^{*}(E)$ is hypercyclically embedded in $G$, then $E$ is hypercyclically embedded in $G$.}

\section{Proofs of Theorems A and  B}

Theorem A is a corollary of the following two general results.

{\bf Theorem 3.1. } {\sl  Let $E$ be a $p$-soluble
 normal subgroup of $G$ and $P$ a Sylow $p$-subgroup
 of $E$. Suppose that $G$ has a complete Hall $\sigma$-set ${\cal H}=\{H_{1}, \ldots , H_{t} \}$  such that
 $H_{1}$ is $p$-supersoluble  of order divisible by $p$. Suppose also that there is a natural number $k$ such that $p^{k} < |P|$ and
 every subgroup of $P$ of order $p^{k}$ and every
cyclic subgroup of $P$ of order 4 (if  $p^k=2$ and $P$ is
 non-abelian) are  $\sigma $-semipermutable in
$G$  with respect to ${\cal H}$. Then $E/O_{p'}(E)$ is hypercyclically embedded in
 $G/O_{p'}(E)$.}

{\bf Theorem 3.2. } {\sl  Let  $P$ be a Sylow $p$-subgroup of $G$. Suppose that
 $G$ has a complete Hall $\sigma$-set ${\cal H}=\{H_{1}, \ldots
, H_{t} \}$  such that $H_{1}$  is  $p$-supersoluble  of order divisible by $p$.
 Suppose also that there is a natural number
$k$ such that $p^{k} < |P|$ and every subgroup of $P$ of order $p^{k}$ and every
 cyclic subgroup of $P$ of order 4 (if $p^k=2$ and $P$ is  non-abelian) are $\sigma $-semipermutable in $G$  with respect to ${\cal H}$.
 Then $G$ is $p$-soluble. }

{\bf Proof of Theorem 3.1.} Assume that this theorem is false and let $G$ be a counterexample
 with $|G|+ |E|$ minimal. Let $|P|=p^{n}$. Then:

(1) {\sl $O_{p'}(N)=1$ for every subnormal subgroup $N$ of $E$.  Hence $O_{p}(G)\ne 1 $.}

Suppose that  for some subnormal subgroup $N$ of $G$ contained in $E$ we have
 $O_{p'}(N)\ne 1$.  Then $O_{p'}(N)$ is subnormal in $G$ and  so
 $O_{p'}(N)\leq O_{p'}(G)$ by Lemma 2.10. On the other hand, by Lemma 2.1(1),  the hypothesis
  holds for   $(G/(E\cap O_{p'}(G)), E/(E\cap O_{p'}(G)))=(G/O_{p'}(E), E/ O_{p'}(E))$.  Hence
 $E/O_{p'}(E)$ is hypercyclically embedded in   $G/O_{p'}(E)$  by the choice of $G$, a
 contradiction. Thus we have (1).

(2) {\sl Let $U=O_{p}(E)$. Then $U$ is not hypercyclically embedded in $G$}.

Assume that $U$ is  hypercyclically embedded in $G$.   Since $E$ is $p$-soluble by
 hypothesis and $O_{p'}(E)=1$ by Claim (1), $U\ne 1$   and $C_{E}(U)\leq U$ by
 the Hall-Higman lemma \cite[Ch.VI, Lemma 6.5]{hupp}. But since $U$ is  hypercyclically
 embedded in $G$, $G/C_{G}(U)$ is strictly $p$-closed  by Lemma 2.4 and so $G/C_G(U)$
 is supersoluble by \cite[Ch.1, Theorem 1.9]{We}. Now in view of the $G$-isomorphism
 $EC_{G}(U)/C_{G}(U)\simeq E/E\cap C_{G}(U)$, we conclude that $E$  is  hypercyclically
 embedded in $G$, a contradiction.

(3) $k > 1$.

Assume that   $k=1$.  We show that in this case  $U$  is  hypercyclically embedded
 in $G$.  Assume that this is false.  Let $U/R$ be a chief factor of $G$. Then by the
 choice of $G$ we have $R$ is hypercyclically embedded in $G$, so for any normal subgroup
 $V$ of $G$ such that $V < U$ we have $V\leq R$ and $U/R$ is not cyclic. Let $B$ be a Thompson
 critical subgroup of $U$ and  $\Omega =\Omega (B)$. We claim that $\Omega  = U$. Indeed,
 if $\Omega < U$, then $ \Omega \leq R$ and so $ \Omega$ is  hypercyclically embedded in $G$.
 Hence $U$   is  hypercyclically embedded in by Lemma 2.5, a contradiction. Thus $\Omega = U$.
 Since $U\leq H_{1}$ and $H_{1}$ is $p$-supersoluble by hypothesis, there is a subgroup
 $L/R\leq U/R$ of order  $p$ such that $L/R$ is normal in $H_{1}/R$. Let $x\in L\setminus R$
 and $H=\langle x \rangle$. Since $\Omega  = U$ and $L\leq U$, $|H|$
is either prime or 4. Then, by hypothesis, $H$ is $\sigma $-semipermutable in $G$  with respect
 to ${\cal H}$. Hence  $HR/R$ is $\sigma $-semipermutable in $G/R$  with respect to
 $\{H_{1}R/R, \ldots , H_{t}R/R \}$ by Lemma 2.1(1). Then,  by Lemma 2.1(5),
 $|G/R:N_{G/R}(HR/R)|=|G/R:N_{G/R}(L/R)|$ is a $\pi (H_{1})$-number. It follows that
  $L/R$ is normal in $G/R$, and so $U/R=L/R$ is cyclic, a contradiction. This shows that $U$
 is hypercyclically embedded in $G$, contrary to Claim (2). Hence we have (3).

(4) {\sl $|N|\leq p^{k}$  for any  minimal normal subgroup $N$ of $G$ contained in $P$} .

Indeed, suppose that $|N| > p^{k}.$ Then  there exists a non-identity proper subgroup $H$
 of $N$ such that $H$ is normal in $H_{1}$ and $H$ is $\sigma$-semipermutable in $G$  with
 respect to ${\cal H}$. But then  $H$ is normal in $G$ by Lemma 2.1(5), which contradicts
 the minimality of $N$.

(5)  {\sl If $P$ is a non-abelian $2$-group, then $k > 2$. }

Assume that $k=2$. We shall show that in this case every subgroup $H$ of $P$ of order 2
 is $\sigma$-semipermutable in   $G$  with respect to ${\cal H}$. This means that $k=1$
 is possible, which will contradicts Claim (3).

First show that for any subgroup $V=A\times  B\leq P$ where $|A|=2=|B|$, if both $V$ and
 $A$ are $\sigma$-semipermutable in $G$  with respect to ${\cal H}$, then $B$ is
$\sigma$-semipermutable in $G$  with respect to ${\cal H}$. Indeed, let $i > 1$ and  $x\in G$.
 Then $AH_{i}^{x}$ and $VH_{i}^{x}$ are  subgroups of $G$ and  $|VH_{i}^{x}:AH_{i}^{x}|=2$.
 Hence $VH_{i}^{x}$ is 2-nilpotent, so $H_{i}^{x}B=H_{i}^{x}B$ since $H_{i}^{x}$ is normal in
 $H_{i}^{x}V$. Similarly, if $V=\langle a \rangle$ is a cyclic
subgroup of order 4, then $\langle a^{2} \rangle$ is $\sigma$-semipermutable in $G$  with
 respect to ${\cal H}$.

Since $P$ is a non-abelian 2-group, $P$ has a cyclic subgroup $H = \langle a \rangle $
 of order 4.  Then $H$ is  $\sigma$-semipermutable in $G$  with respect to ${\cal H}$ by
 hypothesis, so  $A= \langle a^{2} \rangle $ is also  is  $\sigma$-semipermutable in $G$
 with respect to ${\cal H}$. Then every subgroup $B$ of $Z(P)$ of order 2   is
$\sigma$-semipermutable in $G$  with respect to ${\cal H}$, and so every subgroup of $P$ of
 order 2  is  $\sigma$-semipermutable in $G$  with respect to ${\cal H}$.

(6) {\sl If $N$ is a minimal normal subgroup of $G$ contained in $P$, then $(E/N)/O_{p'}(E/N)$
 is hypercyclically embedded in $(G/N)/O_{p'}(E/N)$. }

It is enough to show that the hypothesis holds for $G/N$. Since $E/N$ is $p$-soluble, we can
 assume that $|P/N| > p$.

If either $p > 2$ and $|N|< p^{k}$ or  $p=2$ and $|N| < 2^{k-1}$, then it is clear by Lemma
 2.1(1).  Now let either  $p > 2$ and $|N|= p^{k}$ or $p=2$ and $|N| \in \{2^{k}, 2^{k-1} \} $.

In view of Claim (3),  $k > 1$. Suppose that $|N|=p^{k}$. Then $N$ is non-cyclic and so every
 subgroup of $G$ containing $N$ is not cyclic. Let $N\leq K\leq P$, where $|K:N|=p$. Since $K$
 is non-cyclic, it has a maximal subgroup $L\ne N$. Consider $LN/N$. Since  $L$  is
 $\sigma$-semipermutable in $G$  with respect to ${\cal H}$,  $LN/N$  is also
 $\sigma$-semipermutable in $G/N$ with respect to $\{H_{1}R/R, \ldots , H_{t}R/R \}$ by
 Lemma 2.1(1). Therefore,  if $P/N$ is abelian, the hypothesis is true for $(G/N, P/N)$.
 Next suppose that $P/N$ is a non-abelian $2$-group.

Then $P$ is non-abelian and so $k> 2$ by Claim (5). Since $|P/N|>2$,  $n-k\ge 2.$  We may,
 therefore, let $N\le K\leq V\leq P$ such that $|V:N|=4$, $V/N$ is cyclic and $|V:K|=2$.
 Since $V/N$ is not elementary, $N\nleq \Phi (V)$. Hence for some maximal subgroup $K_{1}$
 of $V$ we have $V=K_{1}N$. Suppose that $K_{1}$ is cyclic. Then
 $|K_1\cap N|=2$ and $2=|V:K_1|=|K_1N:K_1|=|N:K_1\cap N|$. This implies that $|N|=4$. But
 then $k= 2$, a contradiction. Hence $K_{1}$ is not cyclic. Let $S$ and $R$ be two different
 maximal
subgroups of $K_{1}$. Then $K_1=SR.$ If $SN\leq K$ and $RN\leq K$, then $K_{1}=SR\leq K$,
 which contradicts the choice of $K_{1}$. Now since $N/N<K/N<V/N$ where $K/N$ is a maximal
 subgroup of $V/N$, we have that $V/N=K_1N/N=SRN/N=(SN/N)(RN/N).$ But since $V/N$ is cyclic,
 eight $V/N=SN/N$ or $V/N=RN/N$. Without loss of generality, we may assume that $NS=V$.
 Since $S$ is a maximal subgroup of $K_1$ and $K_1$ is a maximal subgroup of $V$,
 $|S|=|N|=p^k$. Then $S$ is $\sigma$-semipermutable in $G$  with respect to ${\cal H}$.
 Hence by Lemma 2.1(1), $V/N$  is  $\sigma$-semipermutable in $G/N$ with respect
 to $\{H_{1}R/R, \ldots , H_{t}R/R \}$. This shows that that the hypothesis is
 true for $(G/N, P/N)$.

Now suppose that  $2^{k-1}=|N|$. If $|N| > 2$, then $N$ is not cyclic and as above one can
 show that every subgroup $\bar{H}$ of $P/N$ with order $2$ and every cyclic subgroup of $P/N$
 of order 4 (if $P/N$ is a non-abelian $2$-group) is $\sigma$-semipermutable in $G/N$   with
 respect to $\{H_{1}R/R, \ldots , H_{t}R/R \}$.
Finally, if  $|N|=2$ and $P/N$ is non-abelian, then $P$ is non-abelian and $k=2$, which
 contradicts Claim (5). Thus (6) holds.

(7) $\Phi (U)=1.$

Assume that for some minimal normal subgroup $N$ of $G$ we have $N\leq \Phi (U)$. Then, by
 Claim (6),  every chief factor of $G/N$ between $ O_{p'}(E/N)$ and $E/N$ is cyclic. Note
 that if $ V/N= O_{p'}(E/N)\ne 1$ and $W$ is a $p$-complement in $V$, then by the Frattini
 argument, $G=VN_{G}(W)=NWN_G(W)=N_{G}(W)$ since $N\leq \Phi (O_{p}(E))\leq \Phi (G)$. Hence
 $W=1$ by Claim (1). Therefore every chief factor of $G$ between $ E$ and $N$ is cyclic.
 Now applying Lemma 2.5, we deduce that $E$  is  hypercyclically embedded in $G$, a
 contradiction. Hence we have (7).

{\sl Final contradiction. } In view of Claims (2) and (7), $U$ is an elementary group and
 for some minimal normal subgroup $N$ of $G$ contained in $U$ we have $|N| > p$. Let $S$ be
 a complement of $N$ in $U$. Since $N\leq H_{1}$ and $|N|\leq p^{k}$ by (4), there are a
 maximal subgroup $V$ of $N$ and a subgroup $W$ of $S$ such that $V$ is normal in $H_{1}$
 and
$|VW|=p^{k}$. Then $VW$ is $\sigma$-semipermutable in $G$  with respect to ${\cal H}$ by
 hypothesis, so $V=VW\cap N$ is normal in $G$ by Lemma 2.1(5). Thus
$V=1$, and so $|N|=p$. This final contradiction completes the proof of the result.

{\bf Proof of Theorem 3.2.} Assume that this theorem is false and let $G$ be a counterexample
 of minimal order.  Without loss of generality we can assume that $P\leq  H_{1}$  and  $H_{i}$
 is a $\sigma _{i}$-group for all $i=1, \ldots , t$. Let $|P|=p^n$ and $V$ be a normal subgroup
 of $G$ such that $G/V$ is a simple group.

(1) {\sl $O_{p'}(N)=1$ for any subnormal subgroup $N$ of $G$} (See Claim (1) in the proof of
 Theorem 3.1).

(2) {\sl $P\nleq N$ for any proper normal subgroup $N$ of $G$} (In view of Lemma 2.1(4), this
 follows from the choice of $G$).

(3) {\sl If the hypothesis holds for $V$,  then $G/V$ is non-abelian, $O_{p}(V)$ is a  Sylow
 $p$-subgroup of $V$ and $O_{p}(V)$  is  hypercyclically embedded in $G$. }

The choice of $G$ implies that $V$ is $p$-soluble. Hence  $V$ is $p$-supersoluble by Theorem A.
 Since $O_{p'}(V)=1$ by Claim (1), $V$ is supersoluble and $O_{p}(V)$ is a Sylow
$p$-subgroup of $V$ by Lemma 2.8. It is clear that $O_{p}(V)$ is normal in $G$, so $O_{p}(V)$
  is  hypercyclically embedded in $G$ by Theorem 3.1.

(4) $k > 1$.

Assume that $k=1$. Then:

(a) {\sl For a Sylow $p$-subgroup $V_{p}$ of $V$   we have $V_{p}\nleq Z_{\infty}(G)$.}

Indeed, assume  that $V_{p}\leq Z_{\infty}(G)$.  By \cite[Ch. IV, Theorem 5.4]{hupp}, $G$ has a $p$-closed
 Schmidt subgroup $A$ and $A=A_{p}\rtimes A_{q}$, where the Sylow  subgroup $A_{p}$ of $A$ is
 of exponent $p$ or exponent 4 (if  $p=2$ and $A_{2}$  is  non-abelian), and if
 $\Phi =\Phi (A_{p})$, then $A_{p}/\Phi$ is a non-central chief factor of $A$. Without loss
 of generality, we may assume that $A_{p}\leq P$. Then
 $V_{p}\cap A\leq Z_{\infty}(A)\cap A_p\leq \Phi$ and so there exists a subgroup $H$
 of $A_{p}$ such that $H\nleq  V$
and $H$ is a cyclic group of order $p$ or of order 4 (if $p=2$ and $A_{2}$  is  non-abelian).
 By hypothesis, $H$ is $\sigma$-semipermutable in $G$, so $HV/V$ is $\sigma$-semipermutable
 subgroup of $G/V$ by Lemma 2.1(1). Note that $G\neq HH_2$ (In fact, if $|H|=p$, it is clear
 since $|P|>p$. If $HH_2=G$ and $H$ is a cyclic group of order 4, then $G$ is $p$-soluble,
 contrary to the choice of $G$).  Hence $G/V$ is not simple by Lemma 2.2, a contradiction.
 Hence we have (a).

(b) {\sl If $|V_{p}|= p$, then $V$ is not $p$-soluble, and so $H_{1}V = G$.}

Indeed, if  $V$ is $p$-soluble, then $V_{p}$ is normal in $G$ by Claim (1). Hence $V_{p} $
 and $C_{G}(V_{p})$ are normal in $G$. Claim (a) implies that $P\leq C_{G}(V_{p}) < G$, which
 contradicts Claim (2). Therefore $V$ is not $p$-soluble. But since the hypothesis holds for
  $H_{1}V$ by Lemma 2.1(2)(3),  the choice of $G$ implies that $H_{1}V = G$.

(c) {\sl  $|V_{p}|\ne p$. Hence the hypothesis holds for $V$ by Lemma 2.1(2) and $|P| > p^{2}$.
 }

Assume that  $|V_{p}|= p$. If $V_{p}=V\cap P\leq \Phi (P)$, then $V$ is $p$-nilpotent by the
 Tate theorem \cite[Ch. IV, Theorem 4.7]{hupp}, contrary to (1). Hence $V_{p}$ has a complement $W$ in $P$.
 Let $L$ be a subgroup of order $p$ in $W$. Assume that $ L < W$. Then the hypothesis holds for
 $VW$ by Lemma 2.1(2)(3), so $VW$ is $p$-soluble, contrary to Claim (b).  Therefore $|W|=p$, so
 $|P|=p^{2}$ and $P=V_pW$ is not cyclic.

Let $E= (H_{2} \cdots H_{t} )^{G}$. Then  in view of Claim (b), we can assume, without of loss
 of generality, that $E\leq V$. We show that  there is a subgroup $W_{0}$ of $P$ order $p$ such
 that $W_{0}\nleq  V$ and  $W_{0}\nleq  C_{G}(E)$. Indeed, suppose that $W\leq  C_{G}(E)$. Note
 that $C_{G}(E)\ne G$ by Claim (1). Hence $V_{p}\nleq C_{G}(E)$ by Claim (2). It follows Claim
 (1) that $C_{G}(E)\cap V=1$. Consequently $G=C_{G}(E)\times V$. Let $W=\langle a\rangle ,$
 $V_{p}=\langle b\rangle $ and $W_{0}=\langle ab\rangle $.  Then
 $W_{0}\cap C_{G}(E)=1= W\cap V$.

Now let $ i > 1$.  Then  $W_{0}H_{i}^{x}=H_{i}^{x}W_{0}$ for all $x\in G$  by hypothesis.
 Let $L=H_{i}^{W_{0}}\cap W_{0}^{H_{i}}$.  Then $L$ is a subnormal subgroup of $G$ by
  \cite[Theorem 7.2.5]{stounh}. Suppose that $L\ne 1$ and let $L_{0}$ be a minimal subnormal subgroup of
 $G$ contained in $L$.
Then $S=L_{0}\cap W_{0}  $ is a Sylow $p$-subgroup of   $L_{0}$  since   $L\leq W_{0}H_{i}$.
 Moreover, in view of Claim (1)  and Lemma 2.10, $S\ne 1$, and so   $W_{0} = S$.
If $L_{0}$ is abelian, then $S= W_{0}  \leq O_{p}(G)$, where $ O_{p}(G)  < P$ by Claim (2).
 Hence   $W_{0} = O_{p}(G)\nleq V$. Consequently $W_{0}\leq  C_{G}(V)\leq C_{G}(E)$. This
  contradiction shows that $L_{0}$ is non-abelian. But then $L_{0}=L_{0}^{G}$  is a minimal
 normal subgroup of $G$ by Claim (2) since $|P|=p^{2}$, which again implies that  $W_{0}\leq
  C_{G}(E)$.
This contradiction shows that  $L=1$. Therefore  for every  $x\in G$ and every $ i > 1$
 we have $(H_{i}^{x})^{W_{0}}\cap W_{0}^{H_{i}^{x}}=1$,
and so  $$[W_{0}, H_{i}^{x}]\leq [(H_{i}^{x})^{W_{0}}, W_{0}^{H_{i}^{x}}] =1.$$
 Therefore  $W_{0}\leq  C_{G}(E)$, a contradiction. Hence  we   have (c).

{\sl Final contradiction for (4). } Let $C=C_{G}(V_{p})$. By Claims (3)  and (c),  $V_{p}$
 is normal in $G$ and it is hypercyclically embedded in $G$. Hence  $G/C$ is
strictly $p$-closed by Lemma 2.4. If  $V_{p}\nleq Z(G)$,  then there is a normal maximal
  subgroup $M$ of $G$ such that $C\leq M$. But since $|P| > p^{2}$, the hypothesis holds for
 $M$, so $M$ is $p$-soluble and so $G$ does.
This contradiction shows that  $V_p\leq Z(G)$, which contradicts Claim (a).  Hence we have (4).

(5) {\sl $|N|\leq p^{k}$  for any  minimal normal subgroup $N$ of $G$ contained in $P$}
 (See Claim (4) in the proof of Theorem 3.1).

(6) $k=n-1$.

Assume that $k < n-1$. Then $VP \ne G$. Indeed, if $VP=G$, then $|G:V|=p$ and the hypothesis
 holds for $V$. Hence $V$ is $p$-soluble by the choice of $G$ and so $G$ is $p$-soluble, a
 contradiction.  By Lemma 2.1(4) the hypothesis holds for $VP$, so $VP$ is $p$-soluble by
 the choice of $G$ since $VP \ne G$. Therefore $V$ is $p$-soluble, so $O_{p}(V)\ne 1$ by
 Claim (1). Let $N$ be a minimal normal subgroup of $G$ contained in $O_{p}(V)$.  It is
 clear that    $N\ne P$. Since $k < n-1$, $|P:N| >p$ by Claim (5). Now repeating  some
 arguments in Claim (6) of the proof of Theorem A one  can show that the hypothesis holds
 for $G/N$, so $G/N$ is $p$-soluble by the choice of $G$. But then $G$ is $p$-soluble, a
 contradiction. Hence we have (6).

(7) {\sl If $O_{p}(G)\ne 1$, then  $P$ is not cyclic. }

Suppose that $P$ is cyclic. Let $L$ be a minimal normal subgroup of $G$ contained in
 $O_{p}(G)\leq P$. Assume that $C_{G}(L)=G$. Then $L\leq Z(G)$. Let $N=N_{G}(P)$.  If
 $P\leq Z(N)$, then $G$ is $p$-nilpotent by Burnside's theorem \cite[Ch. IV, Theorem 2.6]{hupp},
 a contradiction. Hence $N\ne C_{G}(P)$. Let $x\in N\backslash C_{G}(P)$ with $(|x|, |P|)=1$
 and $K=P\rtimes  \langle x \rangle $.  By \cite[Ch. III, Theorem 13.4]{hupp}, $P=[K, P]\times
 (P\cap Z(K))$. Since  $L\leq P\cap Z(K)$ and $P$ is cyclic, it follows that $P= P\cap Z(K)$
 and so $x\in  C_{K}(P)$.  This contradiction shows that $C_{G}(L)\ne G$.

Since $P$ is cyclic,  $|L|=p$. Hence  $G/C_{G}(L)$ is a cyclic group of order dividing $p-1$.
 But then $P\leq C_{G}(L)$, so $C_{G}(L)$ is $p$-soluble by the choice of  $G$. Hence  $G$ is
 $p$-soluble. This contradiction  shows that  we have (7).

(8) {\sl $G\ne PH_{i}$ for any $i>1.$ }

Without lose of generality, assume that $G= PH_{2}$. Let $V_{1}, \ldots , V_{r}$ be the set
 of all maximal subgroups of $P$ and $ D_{i}={V_{i}}^{G}$. Then
$D_{i}={V_{i}}^{PH_{2}}={V_{i}}^{H_{2}}\leq V_{i}H_{2}=H_{2}V_{i}$ by Claim (6).

Suppose that for some $i$, say $i=1$, we have $D_{1}P < G$. Then $D_{1}P$ is $p$-soluble by
the choice of $G$. Hence $O_{p}(G)\ne 1$. By Claim (7), $P$ is not cyclic.
Moreover, for any $ i > 1$, we have that $G=P^{G}=D_{1}D_{i}$. Hence for all such $i>1$, we
 have  that $D_{i}P = G$ and so $D_{i}=V_{i}H_{2}$.
It is also clear  that $V_{2}\cap \cdots \cap V_{r}=\Phi (P)$.
 Let $E= V_{2}H_{2}\cap \cdots \cap V_{r}H_{2}$.
 Then $$P\cap E=(P\cap V_{2}H_{2})  \cap \cdots
\cap  (P\cap V_{r}H_{2}) ==V_{2}(P\cap H_{2})\cap \cdots \cap V_{r}(P\cap H_{2}) =V_{2}\cap
 \cdots \cap V_{r}=\Phi (P).$$  Hence $E$ is $p$-nilpotent by the Tate theorem
\cite[Ch. IV, Thoerem 4.7]{hupp}. It follows that $1 < H_{2}\leq O_{p'}(G)$, contrary to Claim (1). Hence
 we  have (8).

(9) {\sl $P^{G}=G$, so $P\nleq H_{i}^{G} < G$ for all $i > 1$}.

First note that $P^{G}=G$   by Claim (2) and $PH_{i}\ne G$ by Claim (8). If $P$ is not cyclic,
 then $PH^{x}_{i}=H^{x}_{i}P$ for all $x\in G$. Hence $H_{i}^{G} < G$ by Lemma 2.2. Now assume
 that $P$ is cyclic and  $V$ be a maximal subgroup of $P$. Lemma 2.2 implies that either
 $V^{G} < G$ or
$H_{i}^{G} < G$. But if $V^{G} < G$, then $P\nleq V^{G}$ and so $V^{G}\cap P\leq \Phi (P)$.
  Thus   $V^{G}$    is $p$-nilpotent by the Tate theorem
\cite[IV, 4.7]{hupp}, which implies that $V^{G}=V$, contrary to Claim (7). Hence
 $H_{i}^{G} < G$.

{\sl Final contradiction. }  Claim (8) implies that $PH_{i}\ne G$ for all $i=2, \ldots , t$.
 Hence in view of Claim (9), $H_{2}^{G} < G$. Assume that $P=KL$, where $K$ and $L$ are
different maximal subgroups of $P$. Then the hypothesis and claim (6) imply that
  $PH_{i}=KLH_{i}=H_{i}KL=H_{i}P$ for all $i$. On the other hand, the hypothesis holds
 for $PH_{i}$, so $PH_{i}$ is $p$-soluble by the choice of $G$. Now Lemma 2.3 implies that
 $G$ is $p$-soluble. This contradiction shows that $P$ is cyclic. But
$P\nleq H_{2}^{G}$ by Claim (9), so $H_{2}^{G}\cap P\leq \Phi (P)$. Therefore
$H_{2}^{G}$ is  $p$-nilpotent by the Tate theorem \cite[Ch.IV, 4.7]{hupp}. It follows from
 Claim (1) that $H_2^G$ is a  $p$-subgroup. This final
contradiction completes the proof.

{\bf Proof  of Theorem B.} Assume that  this theorem
 is false and let $G$ be a counterexample with $|G|+|E|$
 minimal.

First suppose $X=E$. Let $p$ be the smallest prime dividing
 $|E|$ and $P$ a Sylow $p$-subgroup of $E$. Then  $E$ is
 $p$-nilpotent. Indeed, if $|P|=p$, it follows directly
 from Lemma 2.7. If $|P| >p$, then $E$ is $p$-supersoluble
 by Theorems 3.1 and $3.2$, so $E$ is $p$-nilpotent again by Lemma 2.7.
  Let  $V=O_{p'}(E)$.  Since $V$ is characteristic in $E$, it is normal
 in $G$ and the hypothesis holds for $(G, V)$ and  $(G/V, E/V)$ by Lemma
 2.1(1)(4).

The choice of $G$ and Theorem 3.1 implies that  $P\ne E$. Hence $V\ne 1$,
 so $E/V$  is  hypercyclically embedded in  $G/V$ by the choice of $(G, E)$.
 It is also clear
 that  $V$   is  hypercyclically embedded in $G$. Hence $E$ is  hypercyclically
 embedded in $G$ by the Jordan-H\"older theorem for the chief series, a contradiction.
Therefore in the case, when $X=E$, the theorem is true. Finally, if $X=F^{*}(E)$,
 then the assertion follows from Lemma 2.11. The result is proved.

\section{Applications}

Theorems A,  B,  Theorems 3.1 and 3.2 cover many known results.  Hear we list some
of them.

{\bf Corollary 4.1}  (Gasch\"{u}tz and  N. Ito \cite[Ch. IV, Theorem 5.7]{hupp}). {\sl If
 every minimal subgroup of  $G$  is normal in $G$, then $G$ is soluble and $G'$
 has a normal Sylow 2-subgroup with nilpotent factor group}.

{\bf Proof.} This follows from the fact that $G$ is $p$-supersoluble for all odd
 prime $p$ dividing $|G|$ by Theorem A.

{\bf Corollary 4.2}  (Buckley \cite{Buc}). {\sl If every minimal subgroup
 of a group $G$ of odd order is normal in $G$, then $G$ supersoluble.}

In  view of Example  1.5 we get  from Theorem 3.2 the following results.

{\bf Corollary 4.3}  (Huppert \cite{h2}). {\sl Suppose  that for a Sylow
 $p$-subgroup $P$ of $G$ we have $|P| > p$. Assume that $G$ has a $p$-complement $E$
 such that $E$ permutes with all maximal subgroups of $P$. Then $G$ is $p$-soluble. }

{\bf Corollary 4.4} (Sergienko \cite{Ser}, Borovikov \cite{Bor}) { \sl Suppose  that for a Sylow $p$-subgroup
 $P$ of $G$ we have $|P| > p$.  Assume that $G$ has a $p$-complement $E$  and  there is
 a natural number $k$ such that $p^{k} < |P|$ and every subgroup of $P$ of order $p^{k}$
 permutes with $E$. Suppose also that in the case  when $p=2$ the Sylow 2-subgroups of
 $G$ are abelian. Then $G$ is $p$-supersoluble. }

{\bf Corollary 4.5} (Guo, Shum and Skiba \cite{GuoSS}). {\sl Suppose that
 $G=AT$, where $A$ is a Hall $\pi$-subgroup of $G$ and $T$ a nilpotent supplement of $A$
 in $G$. Suppose that  $A$ permutes with all subgroups of $T$. Then $G$ is $p$-supersoluble
 for each prime $p\not \in \pi$ such that $|T_{p}| > p$ for the Sylow $p$-subgroup $T_{p}$
 of $T$. }

{\bf Proof.} Let $E$ be the Hall $\pi'$-subgroup of $T$. Then every
subgroup $H$  of $E$ permutes with $A^{x}$ for all $x\in G$ by Remark 1.3.
 Hence $H$ is $\sigma$-semipermutable in $G$ with respect to $\{A, E\}$,
so $G$  is $p$-supersoluble by Theorem A.

{\bf Corollary 4.6} (Guo, Shum and Skiba \cite{GuoS163}). {\sl Suppose that
 $G=AT$, where $A$ is a Hall $\pi$-subgroup of $G$ and $T$ a minimal nilpotent supplement
 of $A$
in $G$. Suppose that  $A$ permutes with all maximal subgroups of any Hall  subgroup of $T$.
  Then $G$ is $p$-supersoluble for each prime $p\not \in \pi$ such that $|T_{p}| > p$ for
 the Sylow $p$-subgroup $T_{p}$ of $T$. }

In view of   Example 1.4 we get from Theorem A the following

{\bf Corollary 4.7} (Wei,  Guo \cite{KKK}).  {\sl Let $p$
 be the smallest prime dividing $|G|$ and $P$ be a Sylow $p$-subgroup of $G$.
 If there a subgroup $D$ of $P$ with $1 < |D| < |P|$ such that every subgroup
 $H$ of $P$ with order $|D|$ or order $2|D|$ (if $|D|=2$) is $SS$-quasinormal
 in $G$, then $G$ is $p$-nilpotent. }

From Example 1.4  and Theorem B  we get the following  three results.

 {\bf Corollary   4.8} (Li,  Shen and  Liu \cite{shirong}).
  {\sl Let $\cal F $ be a saturated
formation containing all supersoluble  groups and $ E$  a normal
subgroup
of $G$  such that $G/E \in {\cal F}$.
  Suppose
that for every maximal subgroup of every  non-cyclic Sylow subgroup of $E$ is
 $SS$-quasinormal in $G$. Then  $G \in {\cal F}$.
}

 {\bf Corollary   4.9} (Li,  Shen and   Kong \cite{LiSK}).
  {\sl Let   $ E$ a  normal
subgroup
of $G$  such that $G/E$ is supersoluble.
  Suppose
that for every maximal subgroup of every   Sylow subgroup of $F^{*}(E)$ is
 $SS$-quasinormal in $G$. Then  $G $ is supersoluble.
}

 {\bf Corollary   4.10} (Li,  Shen and   Kong \cite{LiSK}).
  {\sl Let $\cal F $ be a saturated
formation containing all supersoluble  groups and  $ E$ a  normal
subgroup
of $G$  such that $G/E \in {\cal F}$.
  Suppose
that for every maximal subgroup of every   Sylow subgroup of $F^{*}(E)$ is
 $SS$-quasinormal in $G$. Then  $G \in {\cal F}$.
}

\end{document}